\documentclass[reqno]{amsart}
\usepackage{amssymb, amsmath, amsthm, amscd, mathrsfs}

\usepackage{paralist}
\usepackage{cite}
\usepackage{bigints}
\usepackage{dsfont}
\usepackage{empheq}
\usepackage{amssymb}
\usepackage{cases}
\usepackage{enumitem}
\allowdisplaybreaks

\usepackage{caption}
\usepackage{booktabs}

\usepackage{float}
\usepackage[ruled,vlined,linesnumbered,noresetcount]{algorithm2e}

\usepackage{hyperref}

\newtheorem{thm}{Theorem}[section]

\numberwithin{equation}{section}

\newcommand{\ppp}{\partial}

\newcommand{\OOO}{\Omega}

\allowdisplaybreaks

\title[Backward in time semilinear parabolic systems] 
      {Stability of backward-in-time semilinear coupled parabolic systems}

\author{S. E. Chorfi}
\author{M. Yamamoto}

\address{S. E. Chorfi, Cadi Ayyad University, Faculty of Sciences Semlalia, LMDP, UMMISCO (IRD-UPMC), B.P. 2390, Marrakesh, Morocco}
\email{s.chorfi@uca.ac.ma}

\address{M. Yamamoto, Graduate School of Mathematical Sciences, The University of Tokyo, Komaba, Meguro, Tokyo 153-8914, Japan}
\address{Honorary Member of Academy of Romanian Scientists, 
Ilfov, nr. 3, Bucuresti, Romania}
\address{Correspondence member of Accademia Peloritana dei Pericolanti\\
Palazzo Universit\`a, Piazza S. Pugliatti 1 98122 Messina Italy}
\email{myama@ms.u-tokyo.ac.jp}

\makeatletter 
\@namedef{subjclassname@2020}{%
  \textup{2020} Mathematics Subject Classification}
\makeatother

 
\begin{document}
\begin{abstract}
We consider backward problems for semilinear coupled parabolic systems in bounded domains. We prove conditional stability estimates for linear and semilinear systems of strongly coupled parabolic equations involving general semilinearities. The proof of the stability estimates relies on a modified method by Carleman estimates incorporating the simple weight function $e^{\lambda t}$ with a sufficiently large parameter $\lambda$.
\end{abstract}

\maketitle

\section{Introduction}
Let $T>0$ and $N, n\in \mathbb N$. Let $\Omega \subset \mathbb{R}^n$ be a bounded domain with smooth boundary $\partial\Omega$. We set $Q:=\Omega\times (0,T)$ and $\Sigma:=\partial\Omega\times (0,T)$. We consider the semilinear parabolic system in $\mathbf{u}=(u_1, \ldots, u_N)$:
\begin{empheq}[left = \empheqlbrace]{alignat=2}\small
\begin{aligned}
& \partial_t u_\ell(x, t)=\sum_{k=1}^N\sum_{i, j=1}^n \partial_i\left(a_{i j}^{k \ell}(x,t) \partial_j u_k(x, t)\right)+ \sum_{k=1}^N\sum_{i=1}^n b_i^{k \ell}(x, t) \partial_i u_k(x, t) \\
& \hspace{1cm} + \sum_{k=1}^N c^{k \ell}(x, t) u_k(x, t)+ f_\ell(x,t,\mathbf{u}), \hspace{1cm} (x, t) \in Q, \quad 1\le \ell \le N,\\
& \partial_{\nu_A} \mathbf{u}(x, t) +p(x, t) \mathbf{u}(x, t)=0 \text { or } \\
& \mathbf{u}(x, t)=0 \hspace{5.3cm} \text { for }(x, t) \in \Sigma,
\label{eq1}
\end{aligned}
\end{empheq}
where $\nu=(\nu_1,\ldots,\nu_n)$ is the unit outward normal vector field to $\partial\Omega$ and
$$\left[\partial_{\nu_A} \mathbf{u}\right]_{\ell}=\sum_{k=1}^N \sum_{i, j=1}^n a_{i j}^{k \ell}(x,t)\left(\partial_j u_k\right) \nu_i \quad \text { on } \partial \Omega,\quad \ell=1, \ldots, N.$$
Henceforth, for all $1 \leq i, j \leq n, \; 1\le k,\ell \le N$, we assume that
\begin{equation}\label{reg}\small
a_{i j}^{k\ell} \in C^1([0,T];L^{\infty}(\Omega)), \; a_{i j}^{k \ell}=a_{j i}^{k \ell}=a_{i j}^{\ell k}, \; b_i^{k \ell}, c^{k \ell} \in L^\infty(Q), \; p\in W^{1, \infty}\left(0, T ; L^{\infty}(\partial \Omega)\right),
\end{equation}
and that there exists a constant $\sigma>0$ such that
\begin{equation}\label{ella}\small
\sum_{k, \ell=1}^N\sum_{i, j=1}^n a_{i j}^{k \ell}(x,t) \xi_i^k \xi_j^\ell \geq \sigma \sum_{k=1}^N \sum_{i=1}^n |\xi_i^k|^2, \; x \in \overline{\Omega},\;  \xi_i^k \in \mathbb{R}, \; 1 \leq i \leq n, \; 1\le k \le N.
\end{equation}
Semilinear parabolic systems satisfying the symmetry condition in \eqref{reg} arise in fluid dynamics (nonlinear Navier-Stokes equations), see e.g. \cite[Chapter 8]{Ky96}. The ellipticity condition \eqref{ella} is referred to as the strong ellipticity; we refer for instance to \cite{Gi93}.

Let $\beta\in [0,1]$ be fixed. For $\mathbf{f}:=(f_1,\ldots, f_\ell)$, let $\mathbf{F}\colon [0,T]\times H^{\beta}(\Omega)^N \rightarrow L^2(\Omega)^N$, $\mathbf{F}(t,\mathbf{u}(t)):=\mathbf{f}(\cdot,t,\mathbf{u}(\cdot,t))$. We assume that for a.e $(x,t)\in Q$, $\mathbf{f}(x,t;\cdot)$ is a locally Lipschitz continuous function whose Lipschitz constant is independent of $(x,t)$ provided that 
we fix a bounded subset $U$ of $H^{\beta}(\Omega)^N$:
\begin{equation}\label{lip}
    \|\mathbf{f}(x,t,\mathbf{u})-\mathbf{f}(x,t,\mathbf{v})\|_{L^2(\Omega)^N} \le L \|\mathbf{u}-\mathbf{v}\|_{H^{\beta}(\Omega)^N}\quad
\mbox{for all $\mathbf{u}, \mathbf{v}\in U$}.
\end{equation}
Here the constant $L$ is dependent on $U$ but independent of a choice of $\mathbf{u}, \mathbf{v}$.

\textbf{Example.} We consider a one-dimensional version of system \eqref{eq1} with Dirichlet boundary conditions, where $\Omega :=(0,\pi)$ and $f_\ell(x,t,\mathbf{u})=e^{-t}\sin(\partial_x u_\ell(x,t))$ for $1\le \ell\le N$. We can easily check that $\mathbf{f}:=(f_1,\ldots, f_\ell)$ satisfies \eqref{lip} with $\beta =1$.

\subsection*{Backward problem.}
We deal with backward-in-time problems for the semilinear strongly coupled system \eqref{eq1}: for a fixed instant $t_0 \in [0,T)$, determine $\mathbf{u}\left(\cdot, t_0\right)$ from the knowledge of $\mathbf{u}(\cdot, T)$ in $\Omega$.

We set
$$
\mathbf{H}^{2,1}(Q):=\left\{\mathbf{u} \in L^2\left(0,T;L^2(\Omega)^N\right) ; \partial_t \mathbf{u}, \partial_i \mathbf{u}, \partial_i \partial_j \mathbf{u} \in L^2\left(0,T;L^2(\Omega)^N\right) \text{ for } 1 \leq i, j \leq n\right\}.
$$
The main result of conditional stability reads as follows:
\begin{thm}\label{thm1}
We treat two cases:\\
{\bf Case $0<t_0 < T$: (semilinear case)}\\
Let $\mathbf{u}, \mathbf{v}\in \mathbf{H}^{2,1}(Q)$ be two solutions to \eqref{eq1} such that
\begin{equation}\label{ac1}
\hspace{-0.2cm} \exists M>0 \colon \Vert (\mathbf{u}-\mathbf{v})(\cdot,0)\Vert_{H^1(\OOO)^N} \le M, \;\; \|\mathbf{w}(\cdot, t)\|_{H^{\beta}(\Omega)^N} \le M \text{ for } \mathbf{w}=\mathbf{u},\mathbf{v}.
\end{equation}
Then there exist constants $C>0$ and $\theta \in (0,1)$,
dependent on $t_0$ and $M$, such that
\begin{equation} \label{sc1}
\Vert \mathbf{u}(\cdot,t_0)-\mathbf{v}(\cdot,t_0)\Vert_{L^2(\OOO)^N} 
\le C(\Vert \mathbf{u}(\cdot,T)-\mathbf{v}(\cdot,T)\Vert^{\theta}_{H^1(\OOO)^N}
+ \Vert \mathbf{u}(\cdot,T)-\mathbf{v}(\cdot,T)\Vert_{H^1(\OOO)^N}).
\end{equation}

\noindent{\bf Case $t_0=0$: (linear case)}\\
Let $\mathbf{u}, \mathbf{v}\in \mathbf{H}^{2,1}(Q)$ be two solutions to the system \eqref{eq1} with $\mathbf{f}=\mathbf{0}$, such that
$\ppp_t \mathbf{w}, \ppp_t^2 \mathbf{w}\in \mathbf{H}^{2,1}(Q)$ for $\mathbf{w}=\mathbf{u},\mathbf{v}$, and
\begin{equation}\label{ac2}
\exists M=\mathrm{const}>0 \colon \sum_{k=0}^2 \Vert \ppp_t^k (\mathbf{u}-\mathbf{v})(\cdot,0)\Vert_{H^1(\OOO)^N} \le M.
\end{equation}
Then, for any $\alpha \in (0,1)$, there exists a constant $C>0$ such that 
\begin{equation}\label{sc2}
\Vert \mathbf{u}(\cdot,0)-\mathbf{v}(\cdot,0)\Vert_{L^2(\OOO)^N} \le 
C\left( \log \frac{1}{D}\right)^{-\alpha},    
\end{equation}
provided that $\displaystyle D := \sum_{k=0}^2 \Vert \ppp_t^k (\mathbf{u}(\cdot,T)-\mathbf{v}(\cdot,T))\Vert_{H^1(\OOO)^N}$ is sufficiently small.
\end{thm}

\subsection{Literature on stability for nonlinear backward problems}
Conditional stability for backward problems is crucial when dealing with the numerical reconstruction of past or initial data. There are some available methodologies mainly for linear evolution equations; we mention the logarithmic convexity \cite{Pa75}, the weighted energy method \cite{LP61}, and the time analyticity method \cite{KP60}. 
We further refer to a monograph Ames and Straughan \cite{AS}.
However, there are only a few works for nonlinear parabolic equations.

\subsubsection*{Global Lipschitz semilinearities}
First Hölder stability estimates when the semilinear term $F(t,u)$ is Lipschitz continuous in $u$ were obtained by Ames \cite{Am82} using the logarithmic convexity method. In Trong et al. \cite{Tr07}, a one-dimensional backward heat problem for functions $F(x,t,u)$ Lipschitz continuous in $u$ has been regularized by quasi-reversibility and quasi-boundary value methods. Moreover, some error estimates have been given. Trong and Tuan in \cite{Tr08} and \cite{Tu10} gave some error estimates by regularizing the problem via quasi-reversibility or Fourier truncated methods. We also refer to Trong and Tuan \cite{Tr09} for a one-dimensional backward problem regularized by the integral equation method. Nam \cite{Na10} established a Hölder estimate for a regularized solution by using the truncation method.

In their work \cite{HH12}, Hetrick and Hughes established a Hölder continuous dependence for solutions to some approximate well-posed problems. They assume that the semilinearity $F(t,\mathbf{z}(t))$ is globally Lipschitz continuous in both variables. In the paper \cite{Du17}, Duc and Thang have proven a Hölder stability estimate and proposed a stable method.

Note that all the above papers dealt with autonomous equations with positive, self-adjoint operators and global Lipschitz semilinearities.

\subsubsection*{Local Lipschitz semilinearities}
An approximation of a backward nonlinear equation was considered in Long and Ding \cite{Lo94}. The authors used a locally Lipschitz semilinear term $F(u)$. In \cite{Tu14}, Tuan and Trong obtained an error estimate for a regularized solution using the eigenfunction expansion. They further assumed some restrictions on the semilinear term $F(t,u)$ such as $F(t,0)=0$ for all $t$ and some monotnonicity assumption. The paper Trong et al. \cite{Tr15} considered $F(x,t,u)$ locally Lipschitz in $u$ by using the Fourier transform in $\mathbb R^n$. In Hào et al. \cite{Ha18}, the authors used the method of log-convexity to establish Hölder stability estimates for a time-dependent problem with positive self-adjoint 
operators. 

More generally, the above methods are inapplicable to general nonlinear equations with time-dependent coefficients and nonsymmetric operator equations. Our method by Carleman estimates allows us to handle more general nonlinear systems with time-space variable coefficients and the strong coupling. 
We refer to \cite{Y09} for more details on this method. The method requires regularity of the solutions, but this is feasible in our framework thanks to the maximal regularity of parabolic equations, see e.g., \cite{An90}. 
There are some recent works for semilinear parabolic equations that rely on the specific form of the solutions, we refer to Hào et al. \cite{Ha21} 
for a one-dimensional case and to the cited bibliography. Moreover, our method still applies to general nonlinear degenerate parabolic equations. Still, 
we omit the details here and refer to Cannarsa and Yamamoto \cite{Ya23} 
for a degenerate case with semilinear term $F(x,t,u)$. 
As for related work on inverse problems, see Kawamoto \cite{Ka}.
Finally, we refer to \cite{IY14} for coupled parabolic systems where the authors used the uncommon weight function $(t-\eta)^{-2m}$ where $\eta$ is 
a sufficiently small constant and $m\in \mathbb N$ is a large parameter. 

\subsection{Carleman estimate}
For $\mathbf{u}=\left(u_1, \ldots, u_N\right)$ and $\mathbf{v}=\left(v_1, \ldots, v_N\right)$, we set $(\mathbf{u},\mathbf{v})=\displaystyle\sum_{\ell=1}^N u_\ell v_\ell$ and $|\mathbf{u}|^2=\displaystyle\sum_{\ell=1}^N\left|u_\ell\right|^2$. We define an operator $\mathbf{A}(t)$ in $L^2(\Omega)^N$ by $\mathbf{A}(t)=\left(A_1(t), \ldots, A_N(t)\right)$, where 
$$
A_{\ell}(t) \mathbf{u}=\sum_{k=1}^N \sum_{i, j=1}^n \partial_i\left(a_{i j}^{k \ell}(x, t) \partial_j u_k\right)+\sum_{k=1}^N \sum_{i=1}^n b_i^{k \ell}(x, t) \partial_i u_k+\sum_{k=1}^N c^{k \ell}(x, t) u_k
$$
for every $1 \leq \ell \leq N$.

Then \eqref{eq1} can be written as
$$
\begin{cases}\partial_t \mathbf{u}-\mathbf{A}(t) \mathbf{u}=\mathbf{F}(t,\mathbf{u}) & \text { in } Q, \\ \partial_{\nu_A} \mathbf{u} + p \mathbf{u}=0, & \text { on } \Sigma, \\ \mathbf{u}(\cdot, 0)=\mathbf{u}_0 . &
\end{cases}
$$
and the well-posedness can be guaranteed by further assumptions on $\mathbf{A}(t)$, see e.g., \cite{Ta79}. We further set
\begin{equation}
\left\{\begin{array}{l}
\displaystyle (P_{\ell}(x, t, D) \mathbf{u})(x, t)=\partial_t u_{\ell}(x, t)-\sum_{k=1}^N \sum_{i, j=1}^n \partial_i\left(a_{i j}^{k \ell}(x,t) \partial_j u_k(x, t)\right) \\
\displaystyle -\sum_{k=1}^N \sum_{i=1}^n b_i^{k \ell}(x, t) \partial_i u_k(x, t)-\sum_{k=1}^N c^{k \ell}(x, t) u_k(x, t), \\
\displaystyle (P_{\ell}^0(x, t, D) \mathbf{u})(x, t)=\partial_t u_{\ell}(x, t)-\sum_{k=1}^N \sum_{i, j=1}^n \partial_i\left(a_{i j}^{k \ell}(x,t) \partial_j u_k(x, t)\right), \\
x \in \Omega, \quad 0<t<T, \quad 1 \leq \ell \leq N,
\end{array}\right.
\end{equation}
and
$$
\mathbf{P}_0(x, t, D)=\left(P_1^0(x, t, D), \ldots, P_N^0(x, t, D)\right), \mathbf{P}(x, t, D)=\left(P_1(x, t, D), \ldots, P_N(x, t, D)\right).
$$
We consider the weight function
$$\varphi(t)=e^{\lambda t} \ge 1, \qquad t\ge 0,$$
where $\lambda$ is a large parameter suitably chosen. The key ingredient to prove Theorem \ref{thm1} is the following Carleman-type estimate for the strongly coupled parabolic system.

\begin{thm}\label{thm2}
There exists $\lambda_0 >0$ such that for any $\lambda\ge \lambda_0$ we can choose a constant $s_0(\lambda)>0$ so that there exists a constant $C=C(s_0,\lambda_0)>0$ depending on $\Omega, \; \sigma,$ and
$$\small\max_{1 \leq i, j \leq n, 1 \leq k, \ell \leq N}\left\{\left\|a_{i j}^{k \ell}\right\|_{C^1([0,T];L^{\infty}(\Omega))},\left\|b_i^{k \ell}\right\|_{L^{\infty}\left(Q\right)},\left\|c^{k \ell}\right\|_{L^{\infty}\left(Q \right)},\|p\|_{W^{1, \infty}\left(0, T ; L^{\infty}(\partial \Omega)\right)}\right\}$$
such that
\begin{align}
& \int_{Q} \left\{\frac{1}{s \varphi} \left|\partial_t \mathbf{z}\right|^2 +\lambda|\nabla \mathbf{z}|^2+ s \lambda^2 \varphi |\mathbf{z}|^2\right\} e^{2 s \varphi} d x d t \leq C \int_{Q} |\mathbf{P}(x, t, D) \mathbf{z}|^2 e^{2 s \varphi} d x d t \notag\\
& + C \int_{\Omega}\left\{\left(s\lambda \varphi(T) |\mathbf{z}(x, T)|^2 +|\nabla \mathbf{z}(x, T)|^2 \right) e^{2 s \varphi(T)} + \left(s\lambda |\mathbf{z}(x, 0)|^2 + |\nabla \mathbf{z}(x, 0)|^2\right) e^{2 s} \right\} d x \label{car}
\end{align}
for all $s\ge s_0$ and $\mathbf{z} \in H^1\left(0, T ; L^2(\Omega)^N\right) \cap L^2\left(0, T ; H^2(\Omega)^N\right)$ satisfying
\begin{align}
& \partial_{\nu_A} \mathbf{z}(x, t) +p(x, t) \mathbf{z}(x, t)=0 \text { or } \label{bc}\\
& \mathbf{z}(x, t)=0 \quad \text { for any }(x, t) \in \Sigma.
\end{align}
\end{thm}

\section{Proof of the main results} \label{sec2}
\subsection{Proof of Theorem \ref{thm2}}

\begin{proof}
By a density argument, we can assume that $$\mathbf{z} \in C^1\left([0, T] ; C(\overline{\Omega})^N\right) \cap C([0, T]; \left.C^2(\overline{\Omega})^N\right).$$
Here we consider the boundary condition $\partial_{\nu_A} \mathbf{z} + p \mathbf{z}=0$. The case of the Dirichlet boundary condition $\mathbf{z}=0$ on $\Sigma$ is more easy and can be handled similarly. We set
$$
\mathbf{v}=e^{s\varphi} \mathbf{z}, \quad L_{\ell}(x, t, D) \mathbf{v}=e^{s\varphi} P_{\ell}^0(x, t, D)\left(e^{-s\varphi} \mathbf{v}\right) .
$$
Then $\mathbf{v}=\left(v_1, \ldots, v_N\right)$ satisfies the boundary condition \eqref{bc}.
$$
\left(L_{\ell}(x, t, D) \mathbf{v}\right)(x, t)=\partial_t v_{\ell}+\left\{-s\lambda\varphi v_{\ell}-\sum_{k=1}^N \sum_{i, j=1}^n \partial_i\left(a_{i j}^{k \ell}(x,t) \partial_j v_k\right)\right\}, \quad 1 \leq \ell \leq N
$$
and
\begin{align}
\int_{Q} e^{2s\varphi}\left|\mathbf{P}^0(x, t, D) \mathbf{z}\right|^2 d x d t & =  \int_{Q} \sum_{\ell=1}^N\left|L_{\ell}(x, t, D) \mathbf{v}\right|^2 d x d t \notag\\
& \hspace{-1cm} \geq  \int_Q \sum_{\ell=1}^N |\partial_t v_\ell|^2 d x d t - 2 s\lambda \int_{Q} \sum_{\ell=1}^N\left(\partial_t v_{\ell}\right) \varphi v_{\ell} d x d t \notag\\
& \hspace{-1cm} -2 \int_{Q} \sum_{k, \ell=1}^N \sum_{i, j=1}^n\left\{\partial_i\left(a_{i j}^{k \ell}(x,t) \partial_j v_k\right)\right\}\left(\partial_t v_{\ell}\right) d x d t \notag\\
& \hspace{-1cm} =: \int_{Q} |\partial_t \mathbf{v}|^2 d x d t + J_1 + J_2. \label{e1}
\end{align}
In the sequel, $C_j$ will denote generic constants that depend on $$\left\|a_{i j}^{k \ell}\right\|_{C^1([0,T];L^{\infty}(\Omega))},\left\|b_i^{k \ell}\right\|_{L^{\infty}\left(Q\right)},\|c^{k\ell}\|_{L^{\infty}\left(Q\right)},\|p\|_{W^{1, \infty}\left(0, T ; L^{\infty}(\partial \Omega)\right)},\; \Omega, \sigma, N.$$
Hence, through an integration by parts, we have
\begin{align}
J_1=-2 s\lambda & \int_{Q} \varphi \left(\left(\partial_t \mathbf{v}\right), \mathbf{v}\right) d x d t \notag\\
& =-s\lambda \int_{Q} \varphi \left(\partial_t |\mathbf{v}|^2\right) d x d t \notag\\
& =-s\lambda \int_{\Omega}\left[|\mathbf{v}|^2 \varphi\right]_{t=0}^{t=T} d x + s\lambda^2 \int_{Q}|\mathbf{v}|^2 \varphi d x d t \notag\\
& =-s\lambda \int_{\Omega}\left[|\mathbf{z}|^2 e^{2s\varphi} \varphi\right]_{t=0}^{t=T} d x + s\lambda^2 \int_{Q}|\mathbf{z}|^2 e^{2s\varphi} \varphi d x d t . \label{e2}
\end{align}
Using the symmetry of $a_{ij}^{k\ell}$ (see \eqref{reg}), the integration by parts yields
\begin{align}
J_2 &= -  2 \int_{Q} \sum_{k, \ell=1}^N \sum_{i, j=1}^n\left\{\partial_i\left(a_{i j}^{k \ell}(x,t) \partial_j v_k\right)\right\} \partial_t v_{\ell} d x d t \notag\\
& =2 \int_{Q} \sum_{k, \ell=1}^N \sum_{i, j=1}^n a_{i j}^{k \ell}\left(\partial_j v_k\right)\left(\partial_t \partial_i v_{\ell}\right) d x d t-2 \int_{\Sigma} \sum_{\ell=1}^N\left[\partial_{\nu_A} \mathbf{v}\right]_{\ell}\left(\partial_t v_{\ell}\right) d S d t \notag\\
& =2 \int_{Q} \sum_{k, \ell=1}^N \sum_{i, j=1}^n a_{i j}^{k \ell}\left(\partial_j v_k\right)\left(\partial_t \partial_i v_{\ell}\right) d x d t+2 \int_{\Sigma} p(x, t)\left(\mathbf{v}, \partial_t \mathbf{v}\right) d S d t . \label{ee3}
\end{align}
Therefore, by the trace theorem, we have
\begin{align*}
& \left|2 \int_{\Sigma} p(x, t)\left(\mathbf{v}, \partial_t \mathbf{v}\right) d S d t\right|=\left|\int_{\Sigma} p(x, t)\left(\partial_t|\mathbf{v}|^2\right) d S d t\right| \\
& \quad=\left|-\int_{\Sigma}\left(\partial_t p\right) |\mathbf{v}|^2 d S d t+\int_{\partial \Omega}\left[p(x, t)|\mathbf{v}(x, t)|^2\right]_{t=0}^{t=T} d S \right| \\
& \leq C_1\left(\|\mathbf{v}\|_{L^2\left(0, T ; L^2(\partial \Omega)^N\right)}^2+\|\mathbf{v}(\cdot, T)\|_{L^2(\partial \Omega)^N}^2+\|\mathbf{v}(\cdot, 0)\|_{L^2(\partial \Omega)^N}^2\right) \\
& \leq C_2\left(\|\mathbf{v}\|_{L^2\left(0, T ; L^2(\Omega)^N\right)}^2+\|\nabla \mathbf{v}\|_{L^2\left(0, T ; L^2(\Omega)^N\right)}^2+\|\mathbf{v}(\cdot, T)\|_{H^1(\Omega)^N}^2+\|\mathbf{v}(\cdot, 0)\|_{H^1(\Omega)^N}^2\right) .
\end{align*}
Hence we obtain
\begin{align}
& -2 \int_{Q} \sum_{k, \ell=1}^N \sum_{i, j=1}^n\left\{\partial_i\left(a_{i j}^{k \ell}(x,t) \partial_j v_k\right)\right\} \partial_t v_{\ell} d x d t \notag\\
& \geq 2 \int_{Q} \sum_{k, \ell=1}^N \sum_{i, j=1}^n a_{i j}^{k \ell}\left(\partial_j v_k\right)\left(\partial_t \partial_i v_{\ell}\right) d x d t \label{eq25}\\
& \quad-C_2\left(\|\mathbf{v}\|_{L^2\left(0, T ; L^2(\Omega)^N\right)}^2+\|\nabla \mathbf{v}\|_{L^2\left(0, T ; L^2(\Omega)^N\right)}^2+\|\mathbf{v}(\cdot, T)\|_{H^1(\Omega)^N}^2+\|\mathbf{v}(\cdot, 0)\|_{H^1(\Omega)^N}^2\right) . \notag
\end{align}
By \eqref{reg}, we obtain
\begin{align*}
& \sum_{k, \ell=1}^N \sum_{i, j=1}^n a_{i j}^{k \ell}\left(\partial_j v_k\right)\left(\partial_t \partial_i v_{\ell}\right) \\
& =\sum_{k=1}^N \sum_{i=1}^n a_{i i}^{k k}\left(\partial_i v_k\right)\left(\partial_t \partial_i v_k\right)+\sum_{k=1}^N \sum_{i>j} a_{i j}^{k k}\left\{\left(\partial_j v_k\right)\left(\partial_t \partial_i v_k\right)+\left(\partial_i v_k\right)\left(\partial_t \partial_j v_k\right)\right\} \\
& \quad+\sum_{k>\ell} \sum_{i=1}^n a_{i i}^{k \ell}\left\{\left(\partial_i v_k\right)\left(\partial_t \partial_i v_{\ell}\right)+\left(\partial_i v_{\ell}\right)\left(\partial_t \partial_i v_k\right)\right\} \\
& \quad+\sum_{k>\ell} \sum_{i>j}^{k \ell} a_{i j}^{k \ell}\left\{\left(\partial_j v_k\right)\left(\partial_t \partial_i v_{\ell}\right)+\left(\partial_i v_k\right)\left(\partial_t \partial_j v_{\ell}\right)+\left(\partial_j v_{\ell}\right)\left(\partial_t \partial_i v_k\right)+\left(\partial_i v_{\ell}\right)\left(\partial_t \partial_j v_k\right)\right\}\\
&=  \sum_{k=1}^N \sum_{i=1}^n \frac{1}{2} \partial_t\left(\left|\partial_i v_k\right|^2\right) a_{i i}^{k k}+\sum_{k=1}^N \sum_{i>j} a_{i j}^{k k} \partial_t\left(\left(\partial_j v_k\right)\left(\partial_i v_k\right)\right) \\
& +\sum_{k>\ell} \sum_{i=1}^n a_{i i}^{k \ell} \partial_t\left(\left(\partial_i v_k\right)\left(\partial_i v_{\ell}\right)\right)+\sum_{k>\ell} \sum_{i>j} a_{i j}^{k \ell} \partial_t\left(\left(\partial_j v_k\right)\left(\partial_i v_{\ell}\right)+\left(\partial_i v_k\right)\left(\partial_j v_{\ell}\right)\right) \\
&=  \frac{1}{4} \sum_{k, \ell=1}^N \sum_{i, j=1}^n a_{i j}^{k \ell} \partial_t\left(\left(\partial_j v_k\right) \partial_i v_{\ell}+\left(\partial_i v_k\right) \partial_j v_{\ell}\right) .
\end{align*}
Similarly, integration by parts in $t$ yields
\begin{align}
- & 2 \int_{Q} \sum_{k, \ell=1}^N \sum_{i, j=1}^n\left\{\partial_i\left(a_{i j}^{k \ell}(x,t) \partial_j v_k\right)\right\} \partial_t v_{\ell} d x d t \notag\\
& \ge -2 \int_{Q} \sum_{k=1}^N \sum_{i=1}^n \frac{1}{2}\left(\partial_t a_{i i}^{k k}\right)\left|\partial_i v_k\right|^2 d x d t-2 \int_{Q} \sum_{k=1}^N \sum_{i>j}\left(\partial_t a_{i j}^{k k}\right)\left(\partial_j v_k\right)\left(\partial_i v_k\right) d x d t \notag\\
& -2 \int_{Q} \sum_{k>\ell} \sum_{i=1}^n\left(\partial_t a_{i i}^{k \ell}\right)\left(\partial_i v_k\right)\left(\partial_i v_{\ell}\right) d x d t \notag\\
& -2 \int_{Q} \sum_{k>\ell} \sum_{i>j}\left(\partial_t a_{i j}^{k \ell}\right)\left(\left(\partial_j v_k\right)\left(\partial_i v_{\ell}\right)+\left(\partial_i v_k\right)\left(\partial_j v_{\ell}\right)\right) d x d t \notag\\
& +2 \int_{\Omega} \sum_{k=1}^N \sum_{i=1}^n \frac{1}{2}\left[\left|\partial_i v_k\right|^2 a_{i i}^{k k}\right]_{t=0}^{t=T} d x+2 \int_{\Omega} \sum_{k=1}^N \sum_{i>j}\left[a_{i j}^{k k}\left(\partial_j v_k\right)\left(\partial_i v_k\right)\right]_{t=0}^{t=T} d x \notag\\
& +2 \int_{\Omega} \sum_{k>\ell} \sum_{i=1}^n\left[a_{i i}^{k \ell}\left(\partial_i v_k\right)\left(\partial_i v_{\ell}\right)\right]_{t=0}^{t=T}+2 \int_{\Omega} \sum_{k>\ell} \sum_{i>j}\left[a_{i j}^{k \ell}\left(\left(\partial_j v_k\right)\left(\partial_i v_{\ell}\right)+\left(\partial_i v_k\right)\left(\partial_j v_{\ell}\right)\right)\right]_{t=0}^{t=T} d x \notag\\
& -C_2\left(\|\mathbf{v}\|_{L^2\left(Q\right)^N}^2+\|\nabla \mathbf{v}\|_{L^2\left(Q\right)^N}^2+\|\mathbf{v}(\cdot, T)\|_{H^1(\Omega)^N}^2+\|\mathbf{v}(\cdot, 0)\|_{H^1(\Omega)^N}^2\right) \notag\\
&\geq  -C_3 \int_{Q}|\mathbf{z}|^2 e^{2 s\varphi} \varphi d x d t -C_3 \int_{Q}|\nabla \mathbf{z}|^2 e^{2 s\varphi} d x d t \label{e3}\\
& -C_3 \int_{\Omega}\left\{\left(|\nabla \mathbf{z}(x, 0)|^2+|\mathbf{z}(x, 0)|^2\right) e^{2 s}+\left(|\nabla \mathbf{z}(x, T)|^2+|\mathbf{z}(x, T)|^2\right) e^{2 s\varphi(T)}\right\} d x . \notag
\end{align}
By \eqref{ee3}-\eqref{e3} and by enlarging $\lambda$, we obtain
\begin{align}
& \int_{Q} e^{2s\varphi}\left|\mathbf{P}^0(x, t, D) \mathbf{z}\right|^2 d x d t \notag\\
& \ge -s\lambda \int_{\Omega}\left\{|\mathbf{z}(x, T)|^2 e^{2 s\varphi(T)} \varphi(T) - |\mathbf{z}(x, 0)|^2 e^{2s}\right\} d x -C_3 \int_{Q}|\nabla \mathbf{z}|^2 e^{2 s\varphi} d x d t \notag\\
& -C_3 \int_{\Omega}\left\{\left(|\nabla \mathbf{z}(x, T)|^2 +|\mathbf{z}(x, T)|^2\right) e^{2 s\varphi(T)}+\left(|\nabla \mathbf{z}(x, 0)|^2+|\mathbf{z}(x, 0)|^2\right) e^{2 s}\right\} d x \notag\\
& + s\lambda^2 \int_{Q} |\mathbf{z}|^2 e^{2 s\varphi} \varphi d x d t \notag\\
&  \geq  -C_4 s\lambda \int_{\Omega}\left\{|\mathbf{z}(x, T)|^2 e^{2 s\varphi(T)} \varphi(T) + |\mathbf{z}(x, 0)|^2 e^{2s}\right\} d x \notag\\
& -C_4 \int_{\Omega}\left(|\nabla \mathbf{z}(x, T)|^2 e^{2s\varphi(T)} +|\nabla \mathbf{z}(x, 0)|^2 e^{2s}\right) d x-C_4 \int_{Q}|\nabla \mathbf{z}|^2 e^{2s\varphi} d x d t \notag \\
& + s\lambda^2 \int_{Q}|\mathbf{z}|^2 e^{2s\varphi} \varphi d x d t .
\end{align}
Noting that
\begin{align}
& |\mathbf{P}(x, t, D) \mathbf{z}|^2= \sum_{\ell=1}^N\left|P_{\ell}^0(x, t, D) \mathbf{z}-\sum_{k=1}^N \sum_{i=1}^n b_i^{k \ell} \partial_i z_k-\sum_{k=1}^N c^{k \ell} z_k\right|^2 \notag\\
& \quad \geq \sum_{\ell=1}^N \frac{1}{2}\left|P_{\ell}^0(x, t, D) \mathbf{z}\right|^2-C_5\left\{\sum_{\ell=1}^N\left(\left|\sum_{k=1}^N \sum_{i=1}^n b_i^{k \ell} \partial_i z_k\right|^2+\left|\sum_{k=1}^N c^{k \ell} z_k\right|^2\right)\right\} \notag\\
& \quad \geq \frac{1}{2}\left|\mathbf{P}^0(x, t, D) \mathbf{z}\right|^2-\frac{1}{2} C_6|\nabla \mathbf{z}|^2-\frac{1}{2} C_6|\mathbf{z}|^2, \label{eq30}
\end{align}
and taking $s$ and $\lambda$ large enough so that $\frac{s\lambda^2}{2}-\frac{C_6}{2}\ge \frac{s\lambda^2}{4}$, we obtain 
\begin{align}
\int_{Q}e^{2s\varphi}|\mathbf{P}(x, t, D) \mathbf{z}|^2 d x d t &\geq -\frac{C_4 s\lambda}{2} \int_{\Omega}\left\{|\mathbf{z}(x, T)|^2 e^{2s\varphi(T)} \varphi(T) +|\mathbf{z}(x, 0)|^2 e^{2s}\right\} d x \notag\\
& \quad -\frac{C_4}{2} \int_{\Omega}\left\{|\nabla \mathbf{z}(x, T)|^2 e^{2s\varphi(T)}+|\nabla \mathbf{z}(x, 0)|^2 e^{2 s}\right\} d x \notag\\
& \quad -C_7 \int_{Q}|\nabla \mathbf{z}|^2 e^{2s\varphi} d x d t + \frac{s\lambda^2}{4} \int_{Q}|\mathbf{z}|^2 e^{2s\varphi} \varphi d x d t. \notag
\end{align}
Hence,
\begin{align}
s\lambda^2 \int_{Q}|\mathbf{z}|^2 e^{2s\varphi} \varphi d x d t &\le C_8 \int_{Q}e^{2s\varphi}|\mathbf{P}(x, t, D) \mathbf{z}|^2 d x d t + C_8 \int_{Q}|\nabla \mathbf{z}|^2 e^{2s\varphi} d x d t \notag\\
& \quad + C_8 s\lambda \int_{\Omega}\left\{|\mathbf{z}(x, T)|^2 e^{2s\varphi(T)} \varphi(T) +|\mathbf{z}(x, 0)|^2 e^{2s}\right\} d x \notag\\
& \quad + C_8 \int_{\Omega}\left\{|\nabla \mathbf{z}(x, T)|^2 e^{2s\varphi(T)}+|\nabla \mathbf{z}(x, 0)|^2 e^{2 s}\right\} d x .\label{eq28}
\end{align}
Moreover, we have
\begin{align}
&\int_{Q} (\mathbf{P}(x, t, D)\mathbf{z}, \mathbf{z})e^{2s\varphi} d x d t =\int_{Q} \frac{1}{2}\left(\partial_t|\mathbf{z}|^2\right)e^{2s\varphi} d x d t \notag\\
& \quad-\int_0^Te^{2s\varphi}\left(\int_{\Omega} \sum_{k, \ell=1}^N \sum_{i, j=1}^n \partial_i\left(a_{i j}^{k \ell} \partial_j z_k\right) z_{\ell} d x\right) d t \notag\\
& -\int_{Q} \sum_{k, \ell=1}^N \sum_{i=1}^n b_i^{k \ell}\left(\partial_i z_k\right) z_{\ell}e^{2s\varphi} d x d t -\int_{Q} \sum_{k, \ell=1}^N c^{k \ell} z_k z_{\ell}e^{2s\varphi} d x d t .
\end{align}
Let $\gamma>0$ be sufficiently small. By integration by parts, \eqref{ella} and \eqref{bc}, the trace theorem, similar to \eqref{eq25} we obtain
\begin{align}
& \int_{Q}  (\mathbf{P}(x, t, D)\mathbf{z}, \mathbf{z})e^{2s\varphi} d x d t \notag\\
&=  \frac{1}{2} \int_{\Omega}\left(|\mathbf{z}(x, T)|^2e^{2s\varphi(T)}-|\mathbf{z}(x, 0)|^2 e^{2 s}\right) d x -s\lambda \int_{Q} e^{2s\varphi}\varphi |\mathbf{z}|^2 d x d t \notag\\
& +\int_{Q}e^{2s\varphi} \sum_{k, \ell=1}^N \sum_{i, j=1}^n a_{i j}^{k \ell}\left(\partial_j z_k\right)\left(\partial_i z_{\ell}\right) d x d t-\int_{\Sigma}e^{2s\varphi}\left(\partial_{\nu_A} \mathbf{z}, \mathbf{z}\right) d S d t \notag\\
& -\int_{Q}\left(\sum_{k,\ell=1}^N c^{k \ell} z_k z_{\ell}\right)e^{2s\varphi} d x d t-\int_{Q} \sum_{k, \ell=1}^N \sum_{i=1}^n b_i^{k \ell}\left(\partial_i z_k\right) z_{\ell}e^{2s\varphi} d x d t \notag\\
& \geq  -\frac{1}{2} \int_{\Omega}\left(|\mathbf{z}(x, T)|^2e^{2s\varphi(T)}+|\mathbf{z}(x, 0)|^2 e^{2 s}\right) d x \notag\\
& -s\lambda \int_{Q} e^{2s\varphi}\varphi |\mathbf{z}|^2 d x d t+\sigma \int_{Q}|\nabla \mathbf{z}|^2e^{2s\varphi} d x d t \notag\\
& -C_8 \int_0^Te^{2s\varphi}\|\mathbf{z}(\cdot, t)\|^2_{L^2(\partial\Omega)^N} d t \notag\\
& -\int_{Q} \frac{1}{2}e^{2s\varphi} \sum_{k, \ell=1}^N\left\|c^{k \ell}\right\|_{L^{\infty}\left(Q\right)}\left(\left|z_k\right|^2+\left|z_{\ell}\right|^2\right) d x d t \notag\\
& -\int_{Q}e^{2s\varphi} \sum_{k, \ell=1}^N \sum_{i=1}^n\left\|b_i^{k \ell}\right\|_{L^{\infty}\left(Q\right)}\left|\partial_i z_k\right|\left|z_{\ell}\right| d x d t . \label{eq29}
\end{align}
We use the following trace estimate: for any $\varepsilon \in (0,1)$ there exists a constant $C_9(\varepsilon)>0$ such that
\begin{equation}\label{eq210}
    \|\mathbf{z}(\cdot, t)\|_{L^2(\partial\Omega)^N}^2 \leq \varepsilon\|\nabla \mathbf{z}(\cdot, t)\|_{L^2(\Omega)^N}^2+C_9(\varepsilon)\|\mathbf{z}(\cdot, t)\|_{L^2(\Omega)^N}^2, 
\end{equation}
see \cite[Theorem 1.5.1.10]{Gr85}. Moreover
\begin{align}
& \int_{Q}e^{2s\varphi} \sum_{k, \ell=1}^N \sum_{i=1}^n\left\|b_i^{k \ell}\right\|_{L^{\infty}\left(Q\right)}\left|z_{\ell} \| \partial_i z_k\right| d x d t \notag\\
& =\int_{Q}e^{2s\varphi} \sum_{k, \ell=1}^N \sum_{i=1}^n\left(\left\|b_i^{k \ell}\right\|_{L^{\infty}\left(Q\right)}\left|z_{\ell}\right| \sqrt{\frac{N}{\sigma}}\right)\left(\sqrt{\frac{\sigma}{N}}\left|\partial_i z_k\right|\right) d x d t \notag\\
& \leq \int_{Q} \sum_{k, \ell=1}^N \sum_{i=1}^n\left\|b_i^{k \ell}\right\|_{L^{\infty}\left(Q\right)}^2\left|z_{\ell}\right|^2 \frac{N}{2 \sigma}e^{2s\varphi} d x d t +\int_{Q} \sum_{k, \ell=1}^N \sum_{i=1}^n \frac{\sigma}{2 N}\left|\partial_i z_k\right|^2e^{2s\varphi} d x d t\notag\\
& \leq C_{10} \int_{Q}|\mathbf{z}|^2e^{2s\varphi} d x d t+\frac{\sigma}{2} \int_{Q}|\nabla \mathbf{z}|^2e^{2s\varphi} d x d t \notag\\
& \leq C_{10} \int_{Q} |\mathbf{z}|^2 e^{2s\varphi} \varphi  d x d t+\frac{\sigma}{2} \int_{Q}|\nabla \mathbf{z}|^2e^{2s\varphi} d x d t . \label{eq211}
\end{align}
Consequently, \eqref{eq29}-\eqref{eq211} yield
\begin{align}
& \int_{Q} \lambda (\mathbf{P}(x, t, D)\mathbf{z}, \mathbf{z})e^{2s\varphi} d x d t \notag\\
& \quad \geq-\frac{\lambda}{2} \int_{\Omega}\left(|\mathbf{z}(x, T)|^2e^{2s\varphi(T)}+|\mathbf{z}(x, 0)|^2 e^{2 s}\right) d x \notag\\
& \quad-\left(s\lambda^2 +C_{11}\lambda +C_{12} C_9(\varepsilon)\lambda \right) \int_{Q}|\mathbf{z}|^2 e^{2s\varphi} \varphi d x d t \notag\\
& \quad+ \lambda\left(\frac{\sigma}{2}-C_{12} \varepsilon\right) \int_{Q}|\nabla \mathbf{z}|^2e^{2s\varphi} d x d t . \label{eq212}
\end{align}
We choose $\varepsilon \in (0,1)$ such that $\frac{\sigma}{2}-C_{12} \varepsilon=\frac{\sigma}{4}$ ($C_{12}$ can be taken large enough). Since
\begin{align*}
\int_{Q} \lambda (\mathbf{P}(x, t, D)\mathbf{z}, \mathbf{z})e^{2s\varphi} d x d t \leq & \int_{Q} \frac{1}{2}|\mathbf{P}(x, t, D)\mathbf{z}|^2e^{2s\varphi} d x d t +\int_{Q} \frac{\lambda^2}{2}|\mathbf{z}|^2 e^{2s\varphi} \varphi d x d t,
\end{align*}
we obtain, for $s$ and $\lambda$ sufficiently large, 
\begin{align}
\lambda \int_{Q}|\nabla \mathbf{z}|^2e^{2s\varphi} d x d t & \le C_{13}  s\lambda^2\int_{Q} |\mathbf{z}|^2 e^{2s\varphi} \varphi d x d t + C_{13} \int_{Q} |\mathbf{P}(x, t, D)\mathbf{z}|^2e^{2s\varphi} d x d t\notag\\
&  + C_{13} \lambda \int_{\Omega}\left\{|\mathbf{z}(x, T)|^2e^{2s\varphi(T)}+|\mathbf{z}(x, 0)|^2 e^{2 s}\right\} d x. \label{eq3}
\end{align}
Estimating the first term on the right-hand side of \eqref{eq3} in terms of \eqref{eq28}, we obtain
\begin{align}
\lambda \int_{Q}|\nabla \mathbf{z}|^2e^{2s\varphi} d x d t & \le C_{14} \int_{Q} |\mathbf{P}(x, t, D)\mathbf{z}|^2e^{2s\varphi} d x d t + C_{14} \int_{Q}|\nabla \mathbf{z}|^2 e^{2s\varphi} d x d t \notag\\
&  + C_{14} s\lambda \int_{\Omega}\left\{|\mathbf{z}(x, T)|^2 e^{2s\varphi(T)} \varphi(T) +|\mathbf{z}(x, 0)|^2 e^{2s}\right\} d x \notag\\
& + C_{14} \int_{\Omega}\left\{|\nabla \mathbf{z}(x, T)|^2 e^{2s\varphi(T)}+|\nabla \mathbf{z}(x, 0)|^2 e^{2 s}\right\} d x, \label{eq31}
\end{align}
where we employed $\varphi(T)\ge 1$. The second integral on the right-hand side can be absorbed into the left-hand side by taking $\lambda>0$ large enough. Adding up the resulting inequality and \eqref{eq28}, we obtain
\begin{align}
& s\lambda^2 \int_{Q} \varphi |\mathbf{z}|^2e^{2s\varphi} d x d t + \lambda \int_{Q}|\nabla \mathbf{z}|^2e^{2s\varphi} d x d t \notag\\
& \le C_{15} \int_{Q} |\mathbf{P}(x, t, D)\mathbf{z}|^2e^{2s\varphi} d x d t \notag\\
&  + C_{15} s\lambda \int_{\Omega}\left\{|\mathbf{z}(x, T)|^2 e^{2s\varphi(T)} \varphi(T) +|\mathbf{z}(x, 0)|^2 e^{2s}\right\} d x \notag\\
& + C_{15} \int_{\Omega}\left\{|\nabla \mathbf{z}(x, T)|^2 e^{2s\varphi(T)}+|\nabla \mathbf{z}(x, 0)|^2 e^{2 s}\right\} d x .\label{eq32}
\end{align}

In the next step, we will invoke the term in $|\partial_t \mathbf{z}|^2$. Since $\mathbf{z}=e^{-s \varphi} \mathbf{v}$, we have $\partial_t \mathbf{z}=-s \lambda \varphi e^{-s \varphi} \mathbf{v} +e^{-s \varphi} \partial_t \mathbf{v}$, we have
$$
\frac{1}{s \varphi}\left|\partial_t \mathbf{z}\right|^2 e^{2 s \varphi} \leq 2 s \lambda^2 \varphi|\mathbf{v}|^2+\frac{2}{s \varphi}\left|\partial_t \mathbf{v}\right|^2 .
$$
For all large $s>0$ and $\lambda>0$, applying \eqref{e1} yields
\begin{align}
& \int_Q \frac{1}{s \varphi}\left|\partial_t \mathbf{z}\right|^2 e^{2 s \varphi} d x d t \leq \int_Q 2 s \lambda^2 \varphi |\mathbf{v}|^2 d x d t+\int_Q \frac{2}{s \varphi}\left|\partial_t \mathbf{v}\right|^2 d x d t \notag\\
& \leq C_{16} \int_Q s \lambda^2 \varphi|\mathbf{v}|^2 d x d t +C_{16} \int_Q\left|\partial_t \mathbf{v}\right|^2 d x d t \notag\\
& \leq C_{16} \int_Q s \lambda^2 \varphi|\mathbf{v}|^2 d x d t + C_{16} \int_{Q} e^{2s\varphi}\left|\mathbf{P}^0(x, t, D) \mathbf{z}\right|^2 d x d t + C_{16}\left|J_1+J_2\right| . \label{eq33}
\end{align}
Using \eqref{e2} we have
\begin{align*}
\left|J_1\right| \leq s \lambda^2 \int_Q |\mathbf{z}|^2 e^{2s\varphi} \varphi d x d t +s \lambda \int_{\Omega}\left(\varphi(T)|\mathbf{z}(x, T)|^2 e^{2s\varphi(T)}+|\mathbf{z}(x, 0)|^2 e^{2s}\right) d x.
\end{align*}
By \eqref{ee3} and similarly to \eqref{e3}, we obtain 
\begin{align*}
& \quad\left|J_2\right| \le C_{17} \int_{Q}|\mathbf{z}|^2 e^{2 s\varphi} \varphi d x d t +C_{17} \int_{Q}|\nabla \mathbf{z}|^2 e^{2 s\varphi} d x d t \\
& +C_{17} \int_{\Omega}\left\{\left(|\nabla \mathbf{z}(x, 0)|^2+|\mathbf{z}(x, 0)|^2\right) e^{2 s}+\left(|\nabla \mathbf{z}(x, T)|^2+|\mathbf{z}(x, T)|^2\right) e^{2 s\varphi(T)}\right\} d x .
\end{align*}
Substituting these inequalities in \eqref{eq33} and using \eqref{eq32} with \eqref{eq30}, we obtain
\begin{align}
& \int_Q \frac{1}{s \varphi}\left|\partial_t \mathbf{z}\right|^2 e^{2 s \varphi} d x d t \le C_{18} \int_{Q}|\mathbf{z}|^2 e^{2 s\varphi} \varphi d x d t + C_{18} \int_{Q} e^{2s\varphi}\left|\mathbf{P}(x, t, D) \mathbf{z}\right|^2 d x d t \notag\\
& +C_{18} \int_{Q}|\nabla \mathbf{z}|^2 e^{2 s\varphi} d x d t \notag\\
& +  C_{18}\left(s \lambda \varphi(T)\|\mathbf{z}(\cdot, T)\|_{L^2(\Omega)^N}^2+\|\nabla \mathbf{z}(\cdot, T)\|_{L^2(\Omega)^N}^2\right) e^{2 s \varphi(T)} \notag\\
& +  C_{18}\left(s \lambda\|\mathbf{z}(\cdot, 0)\|_{L^2(\Omega)^N}^2+\|\nabla \mathbf{z}(\cdot, 0)\|_{L^2(\Omega)^N}^2\right) e^{2 s} \notag\\
& \leq  C_{19} \int_{Q} e^{2s\varphi}\left|\mathbf{P}(x, t, D) \mathbf{z}\right|^2 d x d t \notag\\
& +C_{19}\left(s \lambda \varphi(T)\|\mathbf{z}(\cdot, T)\|_{L^2(\Omega)^N}^2+\|\nabla \mathbf{z}(\cdot, T)\|_{L^2(\Omega)^N}^2\right) e^{2 s \varphi(T)} \notag\\
& +C_{19}\left(s \lambda\|\mathbf{z}(\cdot, 0)\|_{L^2(\Omega)^N}^2+\|\nabla \mathbf{z}(\cdot, 0)\|_{L^2(\Omega)^N}^2\right) e^{2 s} .
\end{align}
This completes the proof of \eqref{car}.
\end{proof}

\subsection{Proof of Theorem \ref{thm1}}
In the next proofs, $C > 0$ will denote a generic constant that might vary from line to line.

\begin{proof}[Proof of {\bf Case $0<t_0 < T$}]
Setting $\mathbf{z}:=\mathbf{u}-\mathbf{v}$, we obtain the system
\begin{empheq}[left = \empheqlbrace]{alignat=2}\small
\begin{aligned}
& \partial_t z_\ell(x, t)=\sum_{k=1}^N\sum_{i, j=1}^n \partial_i\left(a_{i j}^{k \ell}(x,t) \partial_j z_k(x, t)\right)+ \sum_{k=1}^N\sum_{i=1}^n b_i^{k \ell}(x, t) \partial_i z_k(x, t) \\
& \hspace{0.3cm} + \sum_{k=1}^N c^{k \ell}(x, t) z_k(x, t)+ (f_\ell(x,t,\mathbf{u})-f_\ell(x,t,\mathbf{v})), \hspace{0.5cm} (x, t) \in Q, \quad 1\le \ell \le N,\\
& \partial_{\nu_A} \mathbf{z}(x, t) +p(x, t) \mathbf{z}(x, t)=0 \text { or } \\
& \mathbf{z}(x, t)=0 \hspace{5.3cm} \text { for }(x, t) \in \Sigma.
\label{eq2}
\end{aligned}
\end{empheq}
Applying Theorem \ref{thm2} to \eqref{eq2} yields
\begin{align}\label{es1}
&\bigintssss_Q\left\{\frac{1}{s \varphi} \left|\partial_t \mathbf{z}\right|^2 +\lambda|\nabla \mathbf{z}|^2+ s \lambda^2 \varphi|\mathbf{z}|^2\right\} e^{2 s \varphi} d x d t \leq C \int_Q |\mathbf{f}(x,t,\mathbf{u})-\mathbf{f}(x,t,\mathbf{v})|^2 e^{2 s \varphi} d x d t\notag\\
& + C s \lambda \varphi(T)\|\mathbf{z}(\cdot, T)\|_{H^1(\Omega)^N}^2 e^{2 s \varphi(T)} +C s \lambda\|\mathbf{z}(\cdot, 0)\|_{H^1(\Omega)^N}^2 e^{2 s}
\end{align}
for all large $\lambda$ and $s$. Now we use the assumption on $\mathbf{f}$. By the interpolation inequality, see e.g. \cite[Theorem 1.4.3.3]{Gr85} or \cite[Chapter VII]{Ad75}, there exists a positive constant $C=C(M)$ such that
\begin{align}
\int_Q |\mathbf{f}(x,t,\mathbf{u})-\mathbf{f}(x,t,\mathbf{v})|^2 e^{2 s \varphi} d x d t &\le C\int_0^T \|\mathbf{z}(\cdot, t)\|_{H^{\beta}(\Omega)^N}^2 e^{2s \varphi} dt \notag\\
& \hspace{-0.5cm}\le \int_0^T C\left( \|\mathbf{z}(\cdot, t)\|_{H^1(\Omega)^N}^2 + \|z(\cdot, t)\|_{L^2(\Omega)^N}^2\right) e^{2s \varphi} dt \notag\\
& \hspace{-0.5cm}\le C \int_Q \left(|\nabla \mathbf{z}|^2 + |\mathbf{z}|^2\right) e^{2s \varphi} dx dt.
\end{align}
Therefore, the semilinear term in \eqref{es1} can be absorbed into the left-hand side by taking $\lambda$ large enough. Therefore, by $\varphi(t_0)\le \varphi(t)$, $t\in [t_0,T]$, and using \eqref{ac1}, we obtain
\begin{align*}
& e^{2 s \varphi(t_0)} \bigintssss_{\Omega\times (t_0,T)} \left(\frac{1}{s \varphi} \left|\partial_t \mathbf{z}\right|^2 + s \lambda^2 \varphi|\mathbf{z}|^2\right) d x d t \\
&\le C s \lambda \varphi(T)\|\mathbf{z}(\cdot, T)\|_{H^1(\Omega)^N}^2 e^{2 s \varphi(T)} +C s \lambda M^2 e^{2 s}
\end{align*}
for all large $\lambda$ and $s$. Omitting the dependence on $\lambda=\lambda_0$ large enough and on $T$, we obtain
\begin{align*}
\bigintssss_{\Omega\times (t_0,T)} \left(\frac{1}{s} \left|\partial_t \mathbf{z}\right|^2 + s \varphi|\mathbf{z}|^2\right) d x d t \le C s \|\mathbf{z}(\cdot, T)\|_{H^1(\Omega)^N}^2 e^{2 s (\varphi(T)-\varphi(t_0))} +C s M^2 e^{-2 s \mu(t_0)},
\end{align*}
where $\mu(t):=\varphi(t)-1 \ge 0$. Since $\varphi(t_0)>0$, we deduce
\begin{equation}\label{es2}
\|\partial_t \mathbf{z}\|_{L^2(t_0,T;L^2(\Omega)^N)}^2 \le C s^2 \|\mathbf{z}(\cdot, T)\|_{H^1(\Omega)^N}^2 e^{2 s \varphi(T)} +C s^2 M^2 e^{-2 s \mu(t_0)}
\end{equation}
for all $s$ sufficiently large. Using
$$
\mathbf{z}\left(\cdot, t_0\right)=-\int_{t_0}^T \partial_t \mathbf{z}(\cdot, t) d t + \mathbf{z}(\cdot, T),
$$
there exists a constant $C>0$ such that
\begin{equation}\label{es3}
\left\|\mathbf{z}\left(\cdot, t_0\right)\right\|_{L^2(\Omega)^N}^2 \leq C\left\|\partial_t \mathbf{z}\right\|_{L^2(t_0,T;L^2(\Omega)^N)}^2 +C\|\mathbf{z}(\cdot, T)\|_{L^2(\Omega)^N}^2
\end{equation}
for all $t_0 \in (0, T]$. Since $\varphi(T), \mu(t_0)>0$, $s^2 e^{2 s \varphi(T)} \ge 1$ and $e^{s \varphi(T)} \ge s^2$ for large $s$, combining \eqref{es2} and \eqref{es3}, we obtain
\begin{align*}
& \left\|\mathbf{z}\left(\cdot, t_0\right)\right\|_{L^2(\Omega)^N}^2 \leq C s^2 \|\mathbf{z}(\cdot, T)\|_{H^1(\Omega)^N}^2 e^{2 s \varphi(T)} +C s^2 M^2 e^{-2 s \mu(t_0)} +C \|\mathbf{z}(\cdot, T)\|_{L^2(\Omega)^N}^2 \\
\leq & C s^2 \|\mathbf{z}(\cdot, T)\|_{H^1(\Omega)^N}^2 e^{2 s \varphi(T)} +C s^2 M^2 e^{-2 s \mu\left(t_0\right)} \leq C \|\mathbf{z}(\cdot, T)\|_{H^1(\Omega)^N}^2 e^{3 s \varphi(T)} + C M^2 e^{-s \mu\left(t_0\right)}
\end{align*}
for all large $s \ge s_0$ and all $t_0 \in (0, T]$. We set $K:=C e^{3 s_0 \varphi(T)}$, where $C$ is the same constant in the above inequality. Then
\begin{equation}\label{es4}
\left\|\mathbf{z}\left(\cdot, t_0\right)\right\|_{L^2(\Omega)^N}^2 \leq K \|\mathbf{z}(\cdot, T)\|_{H^1(\Omega)^N}^2 e^{3 s \varphi(T)} +K M^2 e^{-s \mu\left(t_0\right)}    
\end{equation}
for all $s\ge 0$ and all $t_0 \in (0, T]$. By minimizing the right-hand side of \eqref{es4} in $s> 0$ and considering both cases $\|\mathbf{z}(\cdot, T)\|_{H^1(\Omega)^N}\ge M$ and $\|\mathbf{z}(\cdot, T)\|_{H^1(\Omega)^N}\le M$, we finally obtain \eqref{sc1} with $\theta:=\frac{\mu\left(t_0\right)}{3 \varphi(T)+\mu\left(t_0\right)} \in(0,1)$.
\end{proof}

\begin{proof}[Proof of {\bf Case $t_0=0$}]
For simplicity, we assume that the coefficients $a_{ij}^{k \ell}, b_i^{k \ell},  c^{k \ell}$ and $p$, $1 \leq i \leq n$, $1\le k,\ell \le N$, are $t$-independent, otherwise, we can impose further regularity on these coefficients. We will use the previous case $0<t_0<T$. Similarly to previous arguments, we obtain, for $\ppp_t^2 \mathbf{z}$,
\begin{equation}\label{es5}
\|\partial_t^2 \mathbf{z}\|_{L^2(t_0,T;L^2(\Omega)^N)}^2 \le C D^2 e^{3 s \varphi(T)} +C M_1^2 e^{- s \mu(t_0)}
\end{equation}
for all $s\ge 0$ and all $t_0 \in (0,T]$, where $\displaystyle M_1:=\sum_{k=0}^2\left\|\partial_t^k \mathbf{z}(\cdot, 0)\right\|_{H^1(\Omega)^N}.$ By using
$$
\ppp_t \mathbf{z}\left(\cdot, t_0\right)=-\int_{t_0}^T \partial_t^2 \mathbf{z}(\cdot, t) d t + \ppp_t \mathbf{z}(\cdot, T),
$$
and $\left\|\partial_t \mathbf{z}(\cdot, T)\right\|_{L^2(\Omega)^N}^2 \leq C D^2 e^{3 s \varphi(T)}$, we obtain
\begin{align}
\left\|\ppp_t \mathbf{z}\left(\cdot, t_0\right)\right\|_{L^2(\Omega)^N}^2 &\leq C\left\|\partial_t^2 \mathbf{z}\right\|_{L^2(t_0,T;L^2(\Omega)^N)}^2 +C\|\ppp_t \mathbf{z}(\cdot, T)\|_{L^2(\Omega)^N}^2 \notag\\
&\le C D^2 e^{3 s \varphi(T)} +C M_1^2 e^{- s \mu(t_0)}\label{es6}
\end{align}
for all $s\ge 0$ and all $t_0 \in (0,T]$. By means of
$$
\mathbf{z}\left(\cdot, 0\right)=-\int_{0}^T \partial_t \mathbf{z}(\cdot, t) d t + \mathbf{z}(\cdot, T),
$$
and Cauchy-Schwarz inequality, we infer
\begin{align*}
& \|\mathbf{z}(\cdot, 0)\|_{L^2(\Omega)^N}^2 \leq 2 \left\|\int_0^T \partial_t \mathbf{z}\left(\cdot, t\right) d t\right\|_{L^2(\Omega)^N}^2 +2 \|\mathbf{z}(\cdot, T)\|_{L^2(\Omega)^N}^2 \\
\leq & C \int_0^T\left\|\partial_t \mathbf{z}\left(\cdot, t\right)\right\|_{L^2(\Omega)^N}^2 d t +C\|\mathbf{z}(\cdot, T)\|_{L^2(\Omega)^N}^2 .
\end{align*}
Employing \eqref{es6}, we deduce
$$
\|\mathbf{z}(\cdot, 0)\|_{L^2(\Omega)^N}^2 \leq C D^2 \int_0^T e^{3 s \varphi(T)} d t + C M_1^2 \int_0^T e^{- s \mu\left(t\right)} d t + C\|\mathbf{z}(\cdot, T)\|_{L^2(\Omega)^N}^2,
$$
that is,
$$
\|\mathbf{z}(\cdot, 0)\|_{L^2(\Omega)^N}^2 \leq C D^2 e^{C s}+ C M_1^2 \int_0^T e^{- s \mu\left(t\right)} d t.
$$
By the change of variable $\xi:=\mu(t)=e^{\lambda t}-1$, we obtain
$$
\int_0^T e^{- s \mu\left(t\right)} d t=\frac{1}{\lambda} \int_0^{\mu(T)} e^{-2 s \xi} \frac{1}{1+\xi} d \xi \leq \frac{1}{\lambda}\left[-\frac{e^{- s \xi}}{s}\right]_{\xi=0}^{\xi=\mu(T)} =\frac{1-e^{-s\mu(T)}}{\lambda s}\leq \frac{1}{\lambda s}.
$$
Therefore,
$$
\|\mathbf{z}(\cdot, 0)\|_{L^2(\Omega)^N}^2 \leq C D^2 e^{C s}+ \frac{C}{s} M_1^2
$$
for all $s > 0$. Note that we used the same previous argument to extend the inequality for all $s>0$. Without loss of generality, we may assume $D<1$ so that we take $s=\left(\log \frac{1}{D}\right)^\alpha>0$ with arbitrary $\alpha \in (0,1)$. Then
$$
e^{C s} D^2=\exp \left(C\left(\log \frac{1}{D}\right)^\alpha\right) D^2 =\exp \left(C\left(\log \frac{1}{D}\right)^\alpha -2\left(\log \frac{1}{D}\right)\right).
$$
Since $\alpha<1$, there exists a constant $C_1>0$ such that
$$
e^{C\eta^\alpha -2 \eta} \leq C_1 \eta^{-\alpha} \qquad \text { for all } \eta>0,
$$
yielding
$$
e^{C s} D^2 \leq C_1\left(\log \frac{1}{D}\right)^{-\alpha}.
$$
Hence,
$$
C D^2 e^{C_2 s}+\frac{C}{s} M_1^2 \leq C_2\left(\log \frac{1}{D}\right)^{-\alpha}+ C_2 M_1^2\left(\log \frac{1}{D}\right)^{-\alpha}
$$
for some constant $C_2>0$. This completes the proof of Theorem \ref{thm1} for \textbf{Case $t_0=0$}.
\end{proof}

\section*{Acknowledgments}
The second author was supported partly by Grant-in-Aid for 
Scientific Research (A) 20H00117 
and Grant-in-Aid for Challenging Research (Pioneering) 21K18142 of 
Japan Society for the Promotion of Science.


\begin{thebibliography}{99}

\bibitem{Ad75}
R. A. Adams, {\it Sobolev Spaces}, New York: Academic Press, 1975.

\bibitem{Am82}
K. A. Ames, On the comparison of solutions of related properly and improperly posed Cauchy problems for first order operator equations, {\it SIAM J. Math. Anal.}, {\bf 13} (1982), 594-606.

\bibitem{AS}
K. A. Ames and B. Straughan, {\it Non-Standard and Improperly Posed 
Problems}, New York: Academic Press, 1997.

\bibitem{An90}
S. B. Angenent, Nonlinear analytic semiflows, {\it Proc. R. Soc. Edinburgh A}, {\bf 115} (1990) 91-107.

\bibitem{HH12}
B. M. Campbell Hetrick and R. J. Hughes, Continuous dependence on modeling for nonlinear ill-posed problems, {\it J. Math. Anal. Appl.}, {\bf 349} (2009), 420-435.

\bibitem{Ya23}
P. Cannarsa and M. Yamamoto, Stability for backward problems in time for degenerate parabolic equations, (2023), {\it arXiv:} 2305.00525.

\bibitem{Du17}
N. V. Duc and N. V. Thang, Stability results for semi-linear parabolic equations backward in time, {\it Acta. Math. Vietnam.},  {\bf 42} (2017), 99-111.

\bibitem{Gi93}
M. Giaquinta, {\it Introduction to Regularity Theory for Nonlinear Elliptic Systems}, Basel: Birkhäuser, 1993.

\bibitem{Gr85}
P. Grisvard, {\it Elliptic Problems in Nonsmooth Domains}, Monographs and Studies in Mathematics, vol. {\bf 24}, Boston:Pitman, 1985.

\bibitem{Ha21}
D. N. Hào, N. Van Duc and N. T. N. Oanh, Stability results for weak solutions to backward one-dimensional semi-linear parabolic equations with locally Lipschitz source, {\it J. Inverse Ill-Posed Probl.}, {\bf 29} (2021), 499-513.

\bibitem{Ha18}
D. N. Hào, N. V. Duc and N. V. Thang, Backward semi-linear parabolic equations with time-dependent coefficients and local Lipschitz source, {\it Inverse Problems}, {\bf 34} (2018), 055010.

\bibitem{IY14}
O. Y. Imanuvilov and M. Yamamoto, Conditional stability in a backward parabolic system, {\it Appl. Anal.}, {\bf 93} (2014), 2174-2198.

\bibitem{Ka}
A. Kawamoto, Inverse problems for linear degenerate parabolic equations by ``time-like" Carleman estimate, {\it J. Inverse and Ill-Posed Probl.}, {\bf 23} (2015), 1-21.

\bibitem{KP60}
S. G. Krein and O. I. Prozorovskaya, Analytic semigroups and incorrect problems for evolutionary equations, {\it Sov. Math.-Dokl.}, {\bf 133} (1960), 35-38.

\bibitem{Ky96}
P. K. Kythe, {\it Fundamental Solutions for Differential Operators and Applications}, Boston: Birkhäuser, 1996.

\bibitem{LP61}
M. Lees and M. H. Protter, Unique continuation for parabolic differential equations and inequalities, {\it Duke Math. J.}, {\bf 28} (1961), 369-382.

\bibitem{Lo94}
N. T. Long and A. P. N. Ding, Approximation of a parabolic nonlinear evolution equation backwards in time, {\it Inverse Problems}, {\bf 10} (1994), 905-914.

\bibitem{Na10}
P. T. Nam, An approximate solution for nonlinear backward parabolic equations, {\it J. Math. Anal. Appl.}, {\bf 367} (2010), 337-349.

\bibitem{Pa75}
L. E. Payne, {\it Improperly Posed Problems in Partial Differential Equations}, Philadelphia: SIAM, 1975.

\bibitem{Ta79}
H. Tanabe, {\it Equations of Evolution}, London: Pitman; 1979.

\bibitem{Tr15}
D. D. Trong, B. T. Duy and M. N. Minh, Backward heat equations with locally Lipschitz source, {\it Appl. Anal.}, {\bf 94} (2015), 2023-2036.

\bibitem{Tr07}
D. D. Trong, P. H. Quan, T. V. Khanh, N. H. Tuan, A nonlinear case of the 1-D backward heat problem: Regularization and error estimate, {\it Z. fur Anal. ihre Anwend.}, {\bf 26} (2007), 231-245.

\bibitem{Tr08}
D. D. Trong and N. H. Tuan, Stabilized quasi-reversibility method for a class of nonlinear ill-posed problems, {\it Electron. J. Diff. Equ.}, {\bf 84} (2008), 1-12.

\bibitem{Tr09}
D. D. Trong and N. H. Tuan, Regularization and error estimate for the nonlinear backward heat problem using a method of integral equation, {\it Nonlinear Anal.}, {\bf 71} (2009), 4167-4176.

\bibitem{Tu14}
N. H. Tuan and D. D. Trong, On a backward parabolic problem with local Lipschitz source, {\it J. Math. Anal. Appl.}, {\bf 414} (2014), 678-692.

\bibitem{Tu10}
N. H. Tuan and D. D. Trong, A nonlinear parabolic equation backward in time: Regularization with new error estimates, {\it Nonlinear Anal.}, {\bf 73} (2010), 1842-1852.

\bibitem{Y09}
M. Yamamoto, Carleman estimates for parabolic equations and applications, {\it Inverse Problems}, {\bf 25} (2009), 123013.
    
\end{thebibliography}
\end{document}